\documentclass[12pt]{amsart}
\usepackage{amsmath, amssymb}
\usepackage{multirow}

\newtheorem{theorem}{Theorem}[section]
\newtheorem{lemma}[theorem]{Lemma}
\newtheorem{proposition}[theorem]{Proposition}
\newtheorem{corollary}[theorem]{Corollary}
\theoremstyle{remark}
\newtheorem{remark}[theorem]{Remark}

\newcommand{\bQ}{{\mathbb Q}}
\newcommand{\bR}{{\mathbb R}}
\newcommand{\bZ}{{\mathbb Z}}

\def\ord{\operatorname{ord}}
\newcommand{\rad}{{\rm rad}}

\numberwithin{equation}{section}

\begin{document}
\baselineskip=17pt

\title{PRIMITIVE DIVISORS OF CERTAIN ELLIPTIC DIVISIBILITY SEQUENCES}

\author{Paul Voutier}
\address{}
\curraddr{London, UK}
\email{paul.voutier@gmail.com}

\author{Minoru Yabuta}
\address{}
\curraddr{Senri High School, 17-1, 2 chome, Takanodai, Suita, Osaka, 565-0861, Japan}
\email{yabutam@senri.osaka-c.ed.jp, rinri216@msf.biglobe.ne.jp}

\date{}

\dedicatory{}


\begin{abstract}
Let $P$ be a non-torsion point on the elliptic curve $E_{a}: y^{2}=x^{3}+ax$.
We show that if $a$ is fourth-power-free and either $n>2$ is even or $n>1$ is
odd with $x(P)<0$ or $x(P)$ a perfect square, then the $n$-th element of the
elliptic divisibility sequence generated by $P$ always has a primitive divisor.
\end{abstract}

\maketitle

\section{Introduction}

A sequence $C=(C_{n})_{n \ge 1}$ is called a {\it divisibility sequence} if
$C_{m}|C_{n}$ whenever $m|n$. For such a sequence $C$, a prime $p$ is called a
{\it primitive divisor} of the term $C_{n}$ if $p$ divides $C_{n}$ but does not
divide $C_{k}$ for any $0<k<n$. Primitive divisors have been studied by many
authors. In 1892, Zsigmondy~\cite{Zsigmondy} showed that for the sequence
$C_{n}=a^{n}-b^{n}$ the term $C_{n}$ has a primitive divisor for all $n>6$,
where $a$ and $b$ are positive coprime integers. In 1913, Carmichael~\cite{Carm}
showed that if $n>12$ then the $n$-th term of any Lucas sequence has a primitive
divisor in the case of positive discriminant. Ward~\cite{Ward2} and
Durst~\cite{Durst} extended Carmichael's result to Lehmer sequences. In 2001,
Bilu, Hanrot and Voutier~\cite{Bilu} proved that if $n>30$ then every $n$-th
Lucas and Lehmer number has a primitive divisor, and listed all Lucas and Lehmer
numbers without a primitive divisor. The results of Zsigmondy, Carmichael, Ward,
Durst and Bilu, Hanrot and Voutier are all best possible (in the sense that for
$n=6$, $n=12$ and $n=30$, respectively, sequences whose $n$-th element has no
primitive divisor do exist).

Let $E$ be an elliptic curve defined over $\bQ$ and denote by $E(\bQ)$ the
additive group of all rational points on the curve $E$. Let $P \in E(\bQ)$ be a
point of infinite order, and for any non-zero integer $n$ write 
\begin{equation}
\label{eq:110}
x(nP)=\frac{A_{n}(P)}{B_{n}(P)},
\end{equation}
in lowest terms with $A_{n}(P) \in \bZ$ and $B_{n}(P) \in \mathbb{N}$. The
sequence $(B_{n}(P))_{n \geq 1}$ is often called an {\it elliptic divisibility
sequence} in the literature. Strictly speaking, this is not correct as it is
the so-called ``division polynomials'' $\left( \psi_{n}(P) \right)_{n \geq 1}$
(where $x(nP) = \phi_{n}(P)/\psi_{n}^{2}(P)$) which satisfy the required recurrence
relation (see \cite[Section~10.1]{Everest2} and \cite[Exercise~3.7]{Silv2}).

Ward~\cite{Ward1} first studied the arithmetic properties of elliptic divisibility
sequences. Silverman~\cite{Silv4} first showed that for any elliptic curve
$E/\bQ$ in long Weierstrass form and any point $P \in E(\bQ)$ of infinite order,
there exists a positive integer $N_{E,P}$ such that the term $B_{n}(P)$ has a
primitive divisor for all integers $n \ge N_{E,P}$. The bound given by Silverman
is not explicit and not uniform. Everest, Mclaren and Ward~\cite{Everest1}
obtained a uniform and quite small bound beyond which a primitive divisor is
guaranteed for congruent number curves $y^{2}=x^{3}-T^{2}x$ with $T>0$ square-free.

\begin{theorem}[Everest, Mclaren, Ward~\cite{Everest1}]
\label{thm:1.1}
With $E: y^{2}=x^{3}-T^{2}x$ with $T>0$ square-free, let $P \in E(\bQ)$ be a point
of infinite order. If $B_{n}(P)$ does not have a primitive divisor, then\\
\hspace*{3.0mm} {\rm (a)} $n \le 10$ if $n$ is even\\
\hspace*{3.0mm} {\rm (b)} $n \le 3$ if $n$ is odd and $x(P)$ is negative.\\
\hspace*{3.0mm} {\rm (c)} $n \le 21$ if $n$ is odd and $x(P)$ is a rational square.
\end{theorem}

Ingram \cite{Ingram} sharpened the bounds obtained in \cite{Everest1} as follows.

\begin{theorem}[Ingram \cite{Ingram}]
\label{thm:1.2}
Let $E$ and $P$ be as Theorem~$\ref{thm:1.1}$. If $B_{n}(P)$ does not have a
primitive divisor, then $5 \nmid n$, and either $n$ is odd or $n=2$. Furthermore,
if\\
\hspace*{3.0mm} {\rm (a)} $x(P)<0$, or \\
\hspace*{3.0mm} {\rm (b)} $\{x(P),x(P)+T,x(P)-T\}$ contains a rational square,\\ 
then $n \le 2$. 
\end{theorem}

In what follows, $a$ will denote a non-zero integer which is fourth-power-free
and $E_{a}: y^{2}=x^{3}+ax$ will be an elliptic curve. The purpose of this
paper is to generalise the above results on the existence of primitive divisors
to all such $E_{a}$.

This condition on $a$ poses no restriction here, since the minimal model of
any $E_{a}$ is of this form. Furthermore, we require a minimal model, since
otherwise, $B_{n}(P)$ can be without a primitive divisor for arbitrarily large
$n$: for any $n$, let $u=\sqrt{B_{n}(P)}$, $E_{a}'$ be the image of $E_{a}$
under the map $(x,y) \mapsto \left( u^{2}x, u^{3}y \right)$ and $P'$ the image
of $P$ under this map, then $B_{n}(P')=1$.

\begin{theorem}
\label{thm:1.6}
Let $P \in E_{a}(\bQ)$ be a point of infinite order. Let $n$ be a positive
integer and assume that $B_{n}(P)$ does not have a primitive divisor.

If $($i$)$ $n$ is odd and either $x(P)$ is a rational square or $x(P)<0$
or $($ii$)$ $n$ is even, then $n \leq 2$.
\end{theorem} 

\begin{remark}
It is easy to show that there are infinitely many values of $a$ and points
$P \in E_{a}(\bQ)$ such that $x(2P)$ is an integer, so $B_{2}(P)=1$. E.g.,
let $P$ be an integer point on $E_{a}$ for $a \equiv 4 \bmod 16$ by
Lemma~6.1 of \cite{VY1} and the duplication formula. So this theorem is
best-possible.
\end{remark}

\begin{remark}
\label{rem:1.5}
If we remove the conditions on $x(P)$ for odd $n$, then it is again easy to
find infinitely many values of $a$ and points $P \in E_{a}(\bQ)$ such that
$B_{3}(P)=1$. E.g., $P=(30,120) \in E_{-420}(\bQ)$ and
$P=(30,450) \in E_{5850}(\bQ)$.

From our use of Ingram's ideas \cite{Ingram} in Subsection~\ref{subsect:ingram},
it appears that $n \leq 3$ is the best possible bound for $E_{a}$ without any
conditions on $x(P)$.

Indeed, in this paper, we show that if $B_{n}$ does not have a primitive
divisor for any elliptic divisibility sequence generated from any non-torsion
point on any $E_{a}$ for $n>3$, then $5 \nmid n$ and either $n=9, 11, 19$ or
$n>20$.
\end{remark}

%

\begin{remark}
Although $x(2P)$ is a square for $P \in E_{a}(\bQ)$ (see Lemma~6.1(a) in
\cite{VY1}), it is not the case that if $x(P)$ is a perfect square then $P=2Q$
where $Q \in E_{a}(\bQ)$. The significance of this statement is that our result
for $n$ odd and $x(P)$ a perfect square is not a trivial consequence of the
result for even $n$.

\end{remark}

Lastly, we note that Everest, Mclaren and Ward \cite[Theorem~2.4]{Everest1}
proved an analogous result to Theorem~\ref{thm:1.6} (with an absolute, though
unspecified, bound for $n$) for any elliptic curve over $\bQ$ with non-trivial
$2$-torsion, provided that Lang's conjecture holds for the curve. Using Cremona's
elliptic curve data as available via PARI, we found that the elliptic
divisibility sequence defined by
$$
y^{2}+xy = x^{3}-15607620x+23668880400
$$
and $P=(-3780, -167400)$ has no primitive divisor for $n=18$ ($B_{9}=B_{18}=289=17^{2}$)
and satisfies the conditions in Theorem~2.4 of \cite{Everest1}. There are also
five examples up to conductor $130,000$ of such sequences whose $14$-th element
has no primitive divisor.

And although it does not satisfy the conditions in Theorem~2.4 of \cite{Everest1}
the elliptic divisibility sequence defined by
$$
y^{2}+xy = x^{3}-628340x+887106192
$$
and $P=(20824, -3013412)$ has no primitive divisor for $n=21$ ($B_{7}=B_{21}=289=17^{2}$).
We believe this is the largest index, $n$, such that $B_{n}$ has no primitive
divisor for a sequence generated by an elliptic curve with non-trivial $2$-torsion.

Our paper is structured as follows. In the next section, we state our version
of Lang's conjecture for $E_{a}(\bQ)$ as well as bounds on the
difference between the canonical height and the logarithmic (Weil) height for
any point on $E_{a}(\bQ)$. In Section~3, we prove some required results about
elliptic divisibility sequences. Section~4 contains upper bounds on the height
of a point $P \in E_{a}(\bQ)$ generating an elliptic divisibility sequence whose
$n$-th element does not have a primitive divisor, along with their proofs.
Finally Section~5 combines the results of Sections~2--4 to obtain a small
upper bound on such $n$. We then apply some results of Ingram \cite{Ingram}
to complete the proof of Theorem~\ref{thm:1.6}.

\section{Heights of Points on $E_{a}(\bQ)$}

We shall require two results regarding the height of points on elliptic curves.
We require a lower bound in terms of $a$ on the height of nontorsion points on
$E_{a}(\bQ)$. We shall also need bounds on the difference between the canonical
height and the logarithmic (Weil) height for any point on $E_{a}(\bQ)$.

We start by defining these heights.

For a rational point $P \in E(\bQ)$, we define the {\it canonical height} of $P$ by 
$$
\widehat{h}(P)=\frac{1}{2}\lim_{n \to \infty}\frac{h(2^{n}P)}{4^{n}},
$$
where $h(P)=h(x(P))$ is the {\it logarithmic height} of $P$ and
$h(s/t)=\log\max \{|s|, |t|\}$ for a rational number $s/t$ in lowest terms
is the {\it logarithmic height} of $s/t$.

\begin{remark}
This definition of the canonical height follows that in Silverman's book \cite{Silv5}.
This is one-half that found in \cite{Cremona}, as well as one-half that returned
from the height function, ellheight, in PARI (which is important to note here
as we use PARI in some of our calculations).
\end{remark}

\begin{lemma}
\label{lem:lang}
Suppose $a$ is a fourth-power-free integer. Let $P \in E_{a}(\bQ)$ be a
nontorsion point. Then
\begin{equation}
\label{eq:lang2}
\widehat{h}(P) \geq \frac{1}{16}\log |a|
+ \left\{
	\begin{array}{ll}
		 (1/2)\log(2)    & \mbox{if $a>0$ and $a \equiv 1, 5, 7, 9, 13, 15 \bmod 16$} \\
		 (1/4)\log(2)    & \mbox{if $a>0$ and $a \equiv 2, 3, 6, 8, 10, 11, 12, 14 \bmod 16$} \\
		-(1/8)\log(2)    & \mbox{if $a>0$ and $a \equiv 4 \bmod 16$} \\
		 (9/16)\log(2)   & \mbox{if $a<0$ and $a \equiv 1, 5, 7, 9, 13, 15 \bmod 16$} \\
		 (5/16)\log(2)   & \mbox{if $a<0$ and $a \equiv 2, 3, 6, 8, 10, 11, 12, 14 \bmod 16$} \\
		-(1/16)\log(2)   & \mbox{if $a<0$ and $a \equiv 4 \bmod 16$.}
	\end{array}
\right.
\end{equation}
\end{lemma}

\begin{proof}
This is Theorem~1.2 of \cite{VY1}.
\end{proof}

\begin{lemma}
\label{lem:2.6}
For all points $P \in E_{a}(\bQ)$, 
$$
-\frac{1}{4}\log |a|-0.16 \leq \frac{1}{2}h(P)-\widehat{h}(P)
\leq \frac{1}{4}\log |a|+0.26.
$$
\end{lemma}

\begin{proof}
This is Theorem~1.4 of \cite{VY1}.
\end{proof}

\section{Properties of Elliptic Divisibility Sequences}

Let $P \in E_{a}(\bQ)$ be a point of infinite order.
Write 
$$
nP =\left( \frac{A_{n}}{B_{n}}, \; \frac{C_{n}}{B_{n}^{3/2}} \right)
$$
in lowest terms with $A_{n}, C_{n} \in \bZ$ and $B_{n} \in \mathbb{N}$.

\begin{lemma}
\label{lem:2.1}
Let $p$ be any prime divisor of the term $B_{n}$. Then 
$$
\ord_{p}(B_{kn})=\ord_{p}(B_{n})+2\ord_{p}(k).
$$
\end{lemma}

\begin{proof}
This is Lemma~3.1 of \cite{Everest1}.
\end{proof}

\begin{lemma}
\label{lem:2.2}
For any $m, n \in \mathbb{N}$,
$$
\gcd \left( B_{m}, B_{n} \right) = B_{\gcd(m,n)}.
$$
\end{lemma}

\begin{proof}
This is Lemma~3.2 of \cite{Everest1}.
\end{proof}

\begin{lemma} 
\label{lem:2.3}
If the term $B_{n}$ does not have a primitive divisor, then 
\begin{equation}
\label{eq:210}
\log \left( B_{n} \right)
\leq 2 \sum_{p|n}\log(p) +\sum_{p|n} \log\left( B_{n/p} \right).
\end{equation}
\end{lemma}

Here the sums range over prime divisors of $n$.

\begin{proof}
This is the first part of Lemma~3.3 of \cite{Everest1}.
\end{proof}

We let $\omega(n)$ denote the number of distinct prime divisors of $n$.

Further, we define 
\begin{equation}
\label{eq:130}
\rho(n)=\sum_{p|n}p^{-2} \quad \text{and} \quad \eta(n)=2 \sum_{p|n}\log p,
\end{equation}
where the sums range over all prime divisors of $n$ and put
\begin{equation}
\label{eq:K-value}
K=\frac{1}{2}\log |a|+0.52.
\end{equation}

\begin{lemma}
\label{lem:2.7}
If the term $B_{n}$ does not have a primitive divisor, then 
\begin{align}
\log \left( B_{n} \right)
& \leq 2 \sum_{p|n}\log (p) + \sum_{p|n}\left(2\left(\frac{n}{p}\right)^{2}\widehat{h}(P)+K \right) \label{eq:230}\\
& = \eta(n)+2n^{2}\rho(n)\widehat{h}(P)+\omega(n)K. \notag
\end{align}
\end{lemma}

Here the inequality \eqref{eq:230} is analogous to the inequality~(9) of
\cite{Everest1}.

\begin{proof}
Lemma~\ref{lem:2.6} implies that for any prime divisor $p$ of $n$, 
\begin{align}
\log \left( B_{n/p} \right) & \leq h\left(\frac{n}{p}P\right) \notag\\
&\le 2\widehat{h}\left(\frac{n}{p}P\right)+K=2\left(\frac{n}{p}\right)^{2}\widehat{h}(P)+K. \label{eq:220}
\end{align}

The last equality is a property of the canonical height (see Theorem~9.3 of
\cite{Silv2}). Combining \eqref{eq:210} and \eqref{eq:220}, we obtain the lemma.
\end{proof}

\begin{lemma}
\label{lem:2.8}
Let $P \in E_{a}(\bQ)$ be any point of infinite order. Let $m$ and $n$ be
positive integers with $mn$ even.

\noindent
{\rm (a)} If $a \not\equiv 4 \bmod 16$ or $\ord_{2}(x(P)) \neq 1$ or
$\ord_{2}(m)>1$ or $\ord_{2}(n)>1$, then
\begin{equation}
\label{eq:240a}
0 < \left( A_{m}B_{n}-A_{n}B_{m} \right)^{2} = B_{m+n}B_{|m-n|}.
\end{equation}

\noindent
{\rm (b)} If $a \equiv 4 \bmod 16$, $\ord_{2}(x(P))=1$, $\ord_{2}(m)=1$ and
$n$ is odd, then
\begin{equation}
\label{eq:240b}
0 < \left( A_{m}B_{n}-A_{n}B_{m} \right)^{2} = 4B_{m+n}B_{|m-n|}.
\end{equation}
\end{lemma}

\begin{remark}
For $mn$ even, there is one case excluded from this lemma:
$\ord_{2}(m)=\ord_{2}(n)=1$, $a \equiv 4 \bmod 16$ and $\ord_{2}(x(P))=1$.
Here we have
$0 < \left( A_{m}B_{n}-A_{n}B_{m} \right)^{2} = 16B_{m+n}B_{|m-n|}$, but since
we do not need that result in this work, we do not prove that here.

Also the condition that $mn$ is even is a natural one here. Without this
condition, we get additional prime divisors, $p$, of $a$ arising on the right-hand
side when $\ord_{p}(x(P))>0$ from Case~2-b in the proof since $mP$ modulo $p$
is no longer a non-singular point.
\end{remark}

\begin{proof}
Since $P$ is of infinite order, $x(mP) \neq x(nP)$ and hence
$A_{m}B_{n}-A_{n}B_{m} \neq 0$.

For any $P \in E_{a}(\bQ)$ of infinite order and any positive integers
$m$ and $n$, write
$$
mP = \left( x_{m}, y_{m} \right)
   = \left( \frac{A_{m}}{B_{m}}, \frac{C_{m}}{B_{m}^{3/2}} \right), \quad
nP = \left( x_{n}, y_{n} \right)
   = \left( \frac{A_{n}}{B_{n}}, \frac{C_{n}}{B_{n}^{3/2}} \right)
$$
in lowest terms. By the addition formula on the curve $E_{a}$, we have
\begin{eqnarray}
\label{eq:3.8.1}
x(|m \pm n|P) & = & \left( \frac{y_{m} \mp y_{n}}{x_{m}-x_{n}} \right)^{2} - x_{m} - x_{n} \\
\label{eq:3.8.2}
& = & \frac{ \left(C_{m}B_{n}^{3/2} \mp C_{n}B_{m}^{3/2} \right)^{2}}
{B_{m}B_{n} \left( A_{m}B_{n}-A_{n}B_{m} \right)^{2}}
-\frac{A_{m}B_{n}+A_{n}B_{m}}{B_{m}B_{n}}.
\end{eqnarray}

Substituting $y_{m}^{2}=x_{m}^{3}+ax_{m}$ and $y_{n}^{2}=x_{n}^{3}+ax_{n}$ into
(\ref{eq:3.8.1}), we have
\begin{eqnarray*}
&   & x((m+n)P)x(|m-n|P) \\
& = & \frac{\left( \left( x_{m}+x_{n} \right) \left( x_{m}x_{n}+a \right)
-2y_{m}y_{n}\right)\left( (x_{m}+x_{n})(x_{m}x_{n}+a)+2y_{m}y_{n}\right)}
{\left( x_{m}-x_{n}\right)^{4}} \\
& = & \frac{\left(x_{m}+x_{n} \right)^{2} \left( x_{m}x_{n}+a \right)^{2}
-4\left( x_{m}^{3}+ax_{m} \right)\left( x_{n}^{3}+ax_{n} \right)}
{\left( x_{m}-x_{n}\right)^{4}} \\
& = & \frac{\left(x_{m}x_{n}-a \right)^{2}}
{\left( x_{m}-x_{n}\right)^{2}}
= \frac{\left(A_{m}A_{n}-aB_{m}B_{n} \right)^{2}}
{\left( A_{m}B_{n}-A_{n}B_{m}\right)^{2}}.
\end{eqnarray*}

Therefore, we have
\begin{equation}
\label{eq:3.8.2a}
\left( A_{m}B_{n}-A_{n}B_{m}\right)^{2} A_{m+n}A_{|m-n|}
=\left(A_{m}A_{n}-aB_{m}B_{n} \right)^{2}B_{m+n}B_{|m-n|}.
\end{equation}

Put 
\begin{eqnarray*}
G_{m,n} & = & \gcd \left( A_{m+n}A_{|m-n|}, B_{m+n}B_{|m-n|} \right), \\
U_{m,n} & = & \left( A_{m}A_{n} - aB_{m}B_{n} \right)^{2}, \quad
V_{m,n} = \left( A_{m}B_{n} - A_{n}B_{m} \right)^{2}.
\end{eqnarray*}

Then, we have
\begin{equation}
\label{eq:40}
\left( A_{m}B_{n} - A_{n}B_{m} \right)^{2} = \frac{\gcd \left( U_{m,n}, V_{m,n} \right)}{G_{m,n}} B_{m+n}B_{|m-n|}.
\end{equation}

All the above holds without any further conditions on $m$, $n$ and $P$.

To complete the proof, we will prove that $\gcd \left( U_{m,n}, V_{m,n} \right) / G_{m,n}=1$,
under the hypotheses of the lemma.

{\it Step 1: $p=2$.} Under the hypotheses of part~(a) of the lemma, we will
prove that $2 \nmid \gcd \left( U_{m,n}, V_{m,n} \right)$ and under the
hypotheses of part~(b) of the lemma that $\ord_{2} \left( \gcd \left( U_{m,n}, V_{m,n} \right)/G_{m,n} \right)=2$.

First notice that if $B_{m}$ and $B_{n}$ are both even, then $A_{m}A_{n}$ is
odd and so is $A_{m}A_{n}-aB_{m}B_{n}$. Hence $2 \nmid \gcd \left( U_{m,n},V_{m,n} \right)$.

Also if exactly one of $B_{m}$ and $B_{n}$ is even (suppose without loss of
generality that $B_{m}$ is even), then $A_{m}$ is odd and so is
$A_{m}B_{n} - A_{n}B_{m}$. Hence $2 \nmid \gcd \left( U_{m,n},V_{m,n} \right)$.

Therefore, we need only consider the case when both $B_{m}$ and $B_{n}$ are
odd. Furthermore, in this case, if exactly one of $A_{m}$ and $A_{n}$ is even,
then $V_{m,n}$ is odd and again $2 \nmid \gcd \left( U_{m,n},V_{m,n} \right)$.

If $4$ divides either $m$ or $n$ ($m$, without loss of generality), then since
$x(2P)$ is a rational square by Lemma~6.1(a) of \cite{VY1}, it follows from
Lemma~6.1(c) of \cite{VY1} that $B_{4}(P)$, and hence $B_{m}(P)$, is even. We
saw above that this implies that $2 \nmid \gcd \left( U_{m,n},V_{m,n} \right)$.
So we can assume that $\ord_{2}(m), \ord_{2}(n) \leq 1$ and that at least one
of $m$ and $n$ is even (again, suppose $m$ is even).

From Lemma~6.1(c) of \cite{VY1}, unless $a \equiv 4 \bmod 16$ and
$\ord_{2} \left( x(P) \right) = 1$, then
$\ord_{2} \left( B_{m}(P) \right) \geq \ord_{2} \left( B_{2}(P) \right)>0$,
which we saw above implies that $2 \nmid \gcd \left( U_{m,n},V_{m,n} \right)$.
So we can assume that $a \equiv 4 \bmod 16$ and $\ord_{2} \left( x(P) \right)=1$.

We have now handled the hypotheses of part~(a) of the lemma, so we now assume
that $n$ is odd and consider part~(b).

If $A_{m}A_{n}$ is odd, then $U_{m,n}$ is odd too and $2 \nmid \gcd \left( U_{m,n},V_{m,n} \right)$.
So it must be the case that both $A_{m}$ and $A_{n}$ are even.

From the arguments on pages~92--93 of \cite{Silv-Tate}, we can write
$P= \left( b_{1}M^{2}/e^{2}, b_{1}MN/e^{3} \right)$ in lowest terms, where
$a=b_{1}b_{2}$ with $\gcd \left( b_{1}, e \right) = \gcd \left( b_{2}, M \right)
= \gcd \left( e, M \right) = \gcd \left( M, N \right) = \gcd \left( e, N \right)=1$.
Furthermore, since $\ord_{2} \left( x(P) \right)=1$, we have
$2||b_{1}$. Hence $e$ is odd.

Since $2||b_{1}$, we can put $u=2$ in Lemma~6.1(c) of \cite{VY1}. Therefore, we
can write $nP = \left( b_{1}M_{n}^{2}/e_{n}^{2}, b_{1}M_{n}N_{n}/e_{n}^{3} \right)$
in lowest terms, where $e_{n}$ is odd with $\gcd \left( b_{2}, M_{n} \right) =
\gcd \left( M_{n}, N_{n} \right) = 1$. Since $a=b_{1}b_{2} \equiv 4 \bmod 16$
and $2||b_{1}$, we have that $2||b_{2}$. Since $b_{2}$ and $M_{n}$ are coprime,
$M_{n}$ must be odd. Hence $2||A_{n}$.

Since $m$ is even, by Lemma~6.1(a) of \cite{VY1}, $A_{m}$ is a square, so
$4|A_{m}$, since $A_{m}$ is even. Therefore, $4||V_{m,n}$ and so
$\ord_{2} \left( \gcd \left( U_{m,n}, V_{m,n} \right) \right) = 2$,
since $4|U_{m,n}$ too.

Since $B_{4}$ must be even by Lemma~6.1 of \cite{VY1} and $|m \pm n|$ is odd,
it must be the case that both $B_{m+n}$ and $B_{|m-n|}$ are odd (otherwise
$\gcd(4,|m \pm n|)=1$ and so $B_{1}$ is even). Hence $G_{m,n}=1$ and so
if $m \equiv 2 \bmod 4$, $n$ is odd, $a \equiv 4 \bmod 16$ and $\ord_{2}(x(P))=1$,
$\ord_{2} \left( \gcd \left( U_{m,n}, V_{m,n} \right)/ G_{m,n} \right)=2$.

{\it Step 2: $p$, odd.} We will next prove that if $m$ and $n$ are not both
odd, then either $p \nmid \gcd \left( U_{m,n}, V_{m,n} \right)$ or, if 
$p | \gcd \left( U_{m,n}, V_{m,n} \right)$, then
$\ord_{p} \left(\gcd \left( U_{m,n}, V_{m,n} \right) / G_{m,n} \right) = 0$,
for any odd prime $p$.

Put
\begin{eqnarray*}
W_{m,n} & = & \left( C_{m}B_{n}^{3/2} - C_{n}B_{m}^{3/2} \right)^{2}
-\left( A_{m}B_{n} + A_{n}B_{m} \right) \left( A_{m}B_{n} - A_{n}B_{m} \right)^{2}, \\
W_{m,n}' & = & \left( C_{m}B_{n}^{3/2} + C_{n}B_{m}^{3/2} \right)^{2}
- \left( A_{m}B_{n}+A_{n}B_{m} \right) \left( A_{m}B_{n}-A_{n}B_{m} \right)^{2}.
\end{eqnarray*}

Note that $W_{m,n}W_{m,n}'=B_{m}^{2}B_{n}^{2}U_{m,n}V_{m,n}$.

We distinguish four cases.

{\it Case 2-a.}  Assume that $p \nmid W_{m,n}$ and $p \nmid W_{m,n}'$.
Then $p \nmid \gcd \left( U_{m,n}, V_{m,n} \right)$.

{\it Case 2-b.} Assume that $p \mid W_{m,n}$ and $p \mid W_{m,n}'$. We will
prove, by contradiction, that $p \nmid \gcd \left( U_{m,n}, V_{m,n} \right)$.

Suppose that $p$ divides $\gcd \left( U_{m,n}, V_{m,n} \right)$. Then
\begin{eqnarray}
A_{m}A_{n} - aB_{m}B_{n} & \equiv & 0 \bmod p \label{eq:50}, \\
A_{m}B_{n} - A_{n}B_{m}  & \equiv & 0 \bmod p. \label{eq:60}
\end{eqnarray}

If $B_{m} \equiv 0 \bmod p$, then from \eqref{eq:60}, $A_{m}B_{n} \equiv 0 \bmod p$.
Since $A_{m}$ and $B_{m}$ are coprime, we have $A_{m} \not\equiv 0 \bmod p$.
Therefore $B_{n} \equiv 0 \bmod p$. From \eqref{eq:50} we have
$A_{m}A_{n} \equiv 0 \bmod p$, and since $A_{m} \not\equiv 0 \bmod p$, it
follows that $A_{n} \equiv 0 \bmod p$. But this contradicts our assumption
that $A_{n}$ and $B_{n}$ are coprime. Hence $B_{m} \not\equiv 0 \bmod p$.

By the same argument, we obtain $B_{n} \not\equiv 0 \bmod p$. 

Since $mP$ and $nP$ are on $E_{a}$,
\begin{eqnarray}
C_{m}^{2} & \equiv & A_{m}(A_m^{2} + aB_{m}^{2}) \bmod p \label{eq:70} \\
C_{n}^{2} & \equiv & A_{n}(A_n^{2} + aB_{n}^{2}) \bmod p. \label{eq:80}
\end{eqnarray}

From our assumption that $p|W_{m,n}$ and $p|W_{m,n}'$, along with \eqref{eq:60},
we find that
$$
C_{m}B_{n}^{3/2} - C_{n}B_{m}^{3/2} \equiv C_{m}B_{n}^{3/2} + C_{n}B_{m}^{3/2}
\equiv 0 \bmod p,
$$
and hence obtain $2C_{n}B_{m}^{3/2} \equiv 0 \bmod p$. From $B_{m} \not\equiv 0 \bmod p$,
we have $C_{n} \equiv 0 \bmod p$. Furthermore, since $B_{n} \not\equiv 0 \bmod p$,
we also have $C_{m} \equiv 0 \bmod p$. 

First assume that $A_{m}A_{n} \equiv 0 \bmod p$. From \eqref{eq:50}, we have
$aB_mB_n \equiv 0 \bmod p$. Since $B_mB_n \not\equiv 0 \bmod p$, it follows
that $a \equiv 0 \bmod p$. 

Next assume that $A_{m}A_{n} \not\equiv 0 \bmod p$. Then from \eqref{eq:70} and
$C_{m} \equiv 0 \bmod p$, we have 
\begin{eqnarray}
A_{m}^{2} + aB_{m}^{2} \equiv 0 \bmod p. \label{eq:85}
\end{eqnarray}

Multiplying both sides of \eqref{eq:50} by $B_m$ and substituting $A_nB_m \equiv A_mB_n \bmod p$
(from \eqref{eq:60}) and dividing both sides by $B_n$, we obtain $A_m^{2} - aB_{m}^{2} \equiv 0 \bmod p$.
Substracting it from \eqref{eq:85} gives that $2aB_m^2 \equiv 0 \bmod p$ and
therefore $a \equiv 0 \bmod p$.

In both cases $a \equiv 0 \bmod p$ and so from Proposition~VII.5.1 of \cite{Silv2},
$E_{a}$ has bad reduction at $p$ (additive reduction, in fact) and the reduction
of $mP$ modulo $p$ is singular.

On the other hand, from Lemma~4.1(b) of \cite{VY1}, $mP$ modulo
$p$ is non-singular, since $m$ is even. This is a contradiction. Hence $p$ does
not divide $\gcd \left( U_{m,n}, V_{m,n} \right)$. 

{\it Case 2-c.} Assume that $p \mid W_{m,n}$ and $p \nmid W_{m,n}'$. Assume
that $p$ divides $\gcd \left( U_{m,n}, V_{m,n} \right)$ and let
$\ord_{p} \left( \gcd \left( U_{m,n}, V_{m,n} \right) \right) =\alpha>0$.

First assume that $\ord_{p} \left( V_{m,n} \right) \geq \ord_{p} \left( U_{m,n} \right)$.
Then we can put $\ord_{p}\left( U_{m,n} \right)=\alpha$ and $\ord_{p}\left( V_{m,n} \right)=\alpha+\beta$
with $\beta \geq 0$. From \eqref{eq:3.8.2} we obtain
\begin{equation}
\label{eq:3.90}
A_{|m-n|}B_{m}B_{n}V_{m,n} = B_{|m-n|}W_{m,n}'.
\end{equation}

Since $A_{|m-n|}$ and $B_{|m-n|}$ are coprime and $p$ does not divide $W_{m,n}'$,
we have that $\ord_{p} \left( A_{|m-n|} \right)=0$. We saw, in the case~2, that
if $p$ divides $\gcd \left( U_{m,n}, V_{m,n} \right)$, then both $B_{m}$ and $B_{n}$ are prime to $p$.
Hence, from \eqref{eq:3.90} we obtain $\ord_{p} \left( B_{|m-n|} \right)=\alpha+\beta$. 

From \eqref{eq:3.8.2a}, we can write
$$
A_{m+n}A_{|m-n|}V_{m,n} = B_{m+n}B_{|m-n|}U_{m,n}.
$$

Here $\ord_{p}\left( U_{m,n} \right)=\alpha$, $\ord_{p}\left( V_{m,n} \right)=\alpha+\beta$, $\ord_{p} \left( A_{|m-n|} \right)=0$
and $\ord_{p} \left( B_{|m-n|} \right)=\alpha+\beta$. Since $A_{m+n}$ and $B_{m+n}$
are coprime, we obtain that $\ord_{p} \left( A_{m+n} \right)=\alpha$ and
$\ord_{p} \left( B_{m+n} \right)=0$. Hence
$$
\ord_{p} \left( G_{m,n} \right)
=\ord_{p}\left( \gcd \left( A_{m+n}A_{|m-n|}, B_{m+n}B_{|m-n|} \right) \right) =\alpha
$$
and $\ord_{p} \left( \gcd \left( U_{m,n}, V_{m,n} \right) / G_{m,n} \right)=0$.

For the case $\ord_{p}\left( U_{m,n} \right) \geq \ord_{p}\left( V_{m,n} \right)$,
in the same way, we again obtain $\ord_{p}\left( G_{m,n} \right) = \alpha$.

It follows that $\ord_{p}(\gcd \left( U_{m,n}, V_{m,n} \right) / G_{m,n}) = 0$. 

{\it Case 2-d.} Assume that $p \nmid W_{m,n}$ and $p \mid W_{m,n}'$. Then by the
same argument as in Case~2-c, if $p | \gcd \left( U_{m,n}, V_{m,n} \right)$, then
we have $\ord_{p} \left( \gcd \left( U_{m,n}, V_{m,n} \right) / G_{m,n} \right) = 0$.

{\it Step 3: $G_{m,n}$.} We show for any prime, $p$, if
$p \nmid \gcd \left( U_{m,n}, V_{m,n} \right)$, then $p \nmid G_{m,n}$.

Put $\gcd \left( B_{m}, B_{n} \right) = g_{m,n}$. Then we can write $B_{m}=g_{m,n}b_{m}$
and $B_{n}=g_{m,n}b_{n}$, where $b_{m}$ and $b_{n}$ are coprime. Since
$g_{m,n}^{3}|W_{m,n}$ (recall the definition of $W_{m,n}$ from Step~2 as the
common numerator of the expression in \eqref{eq:3.8.2}), we can write
$W_{m,n}=g_{m,n}^{3}W_{m,n,1}$, and from \eqref{eq:3.8.2}, we have
$$
x((m+n)P) = \frac{W_{m,n,1}}{g_{m,n}b_{m}b_{n} \left( A_{m}b_{n}-A_{n}b_{m} \right)^{2}}.
$$

From Lemma~\ref{lem:2.2}, we have
$\gcd \left( B_{m+n}, B_{m} \right) = B_{\gcd(m+n,m)} = B_{\gcd(m,n)}=g_{m,n}$
and similarly, $\gcd \left( B_{m+n}, B_{n} \right) =g_{m,n}$.

Therefore, we can write $B_{m+n}=g_{m,n}b_{m+n}$ where $b_{m+n}$ is prime to
both $b_{m}$ and $b_{n}$. Therefore neither of $b_{m}$ and $b_{n}$ divide
$B_{m+n}$. Hence
\begin{equation}
\label{eq:yabuta-1}
B_{m+n} | \left( A_{m}B_{n} - A_{n}B_{m} \right)^{2}.
\end{equation}

Now assume that $p \nmid \gcd \left( U_{m,n}, V_{m,n} \right)$. We will show by
contradiction that $p$ does not divide $G_{m,n}$.

Recall that $G_{m,n} = \gcd \left( A_{m+n}A_{|m-n|}, B_{m+n}B_{|m-n|} \right)$.
Suppose that $p$ divides $G_{m,n}$. Since $\gcd \left( A_{m+n}, B_{m+n} \right)
= \gcd \left( A_{|m-n|}, B_{|m-n|} \right) = 1$, without loss of generality, we
may assume that $p$ divides  both $A_{|m-n|}$ and $B_{m+n}$.

Let $\ord_{p} \left( A_{|m-n|} \right) = s >0$, $\ord_{p} \left( B_{m+n} \right) = t >0$
and $2\ord_{p} \left( A_{m}B_{n}-A_{n}B_{m} \right) = u >0$. From \eqref{eq:yabuta-1},
we have $t \leq u$. On the other hand, from \eqref{eq:3.8.2a} we have
$u+s=t$ since $p \nmid \gcd \left( U_{m,n}, V_{m,n} \right)$. Hence $s \leq 0$.
This contradiction allows us to conclude that $p$ does not divide $G_{m,n}$.

From these three steps, it follows that $0 < \left( A_{m}B_{n} - A_{n}B_{m} \right)^{2}
= B_{m+n}B_{|m-n|}$, as desired.
\end{proof}

\section{Upper Bounds}

Recall our definitions $\rho(n)$ and $\eta(n)$ from \eqref{eq:130},
as well as $K$ from (\ref{eq:K-value}), and set
\begin{equation}
\label{eq:L-value}
L=\frac{1}{2}\log |a|+0.32.
\end{equation}

We obtain the following proposition. 

\begin{proposition}
\label{prop:1.3}
Let $P \in E_{a}(\bQ)$ be a point of infinite order.
Let $n$ be a positive integer and assume that the term $B_{n}(P)$ does not have
a primitive divisor. \\
{\rm (a)} If $n$ is odd and $x(P)$ is a rational square or if $n$ is even,
write $n=2^{e}N$ where $e$ is a non-negative integer and $N$ is an odd integer,
then $n=1,2,4$ or $N \geq 3$ and
$$
0 < 2\left( \frac{1}{3}-\frac{1}{3N^{2}}-\rho(n) \right) \widehat{h}(P)n^{2}
\leq \eta(n) + \omega(n)K + K + L.
$$

\noindent
{\rm (b)} Let $p$ be an odd prime. If $n$ is odd, divisible by $p$ and $x(P)$
is a rational square, or if $n$ is even and divisible by $p$, then 
$$
0 < 2\left( \frac{(p+1)^{2}}{4p^{2}}-\rho(n) \right) \widehat{h}(P)n^{2}
\leq \eta(n) + \omega(n)K + L.
$$
\end{proposition}

\begin{remark}
Note that part~(b) can be improved for $n$ even and $a<0$:
$$
0 < 2 \left( \frac{5p^{2}+6p+5}{16p^{2}} - \rho(n) \right) \widehat{h}(P)n^{2}
\leq \eta(n) + \omega(n)K + 2L + \log |a| + 0.385.
$$

However, due to the restriction to $a<0$, we will not use it here and so we do
not include the proof (it is available on request). We mention it so other
researchers know that improvements can be made.
\end{remark}

By using estimates for $\rho(n)$, $\omega(n)$ and $\eta(n)$, we obtain the
following corollary.

\begin{corollary}
\label{cor:1.4}
Let $P \in E_{a}(\bQ)$ be a point of infinite order.\\
{\rm (a)} Let $n \geq 3$ be an odd integer and assume that $x(P)$ is a rational
square. If $B_{n}(P)$ does not have a primitive divisor, then 
$$
0.484 \widehat{h}(P)n^{2}
< 2 \log (n) + \frac{1.3841 \log (n)}{\log \log (n)}K + K + L.
$$

{\rm (b)} Let $n$ be a positive even integer not divisible by $5$. If $B_{n}(P)$
does not have a primitive divisor, then either $n \leq 4$ or $n$ is not a power
of $2$ and
$$
0.049 \hat{h}(P)n^{2}
< 2 \log (n) + \frac{1.3841 \log (n)}{\log \log (n)}K + L.
$$
\end{corollary}

\subsection{Proof of Proposition~\ref{prop:1.3}}

We are now ready to prove Proposition~\ref{prop:1.3}. Our proof is based upon
ideas found in \cite{Everest1}.

\subsubsection{Proof of part~(a)}

Assume that either $n>1$ is an odd integer and $x(P)$ is a rational square
or $n$ is even.

If $B_{2^{m}}(P)$ does not have a primitive divisor, then $m \leq 2$
(see Theorem~1.2 of \cite{Yabuta}). Hence we may assume that $n$ is not a power
of two, and write $n=2^{e}N$, where $e$ is a non-negative integer and $N$ is an
odd integer with $N \geq 3$. 

Write $N=3k+r$ with $r=0, \pm 1$, and put $m=2^{e}(2k+r)$ and $m'=2^{e}k$.
Since $N>1$, we have $k>0$ and so $m'>0$ and $m-m'=2^{e}(k+r)>0$.
Also $n=m+m'$.

If $r=\pm 1$, then $k$ is even and $2k+r$ is odd. If $n$ is odd, then $m$ is
odd and $m'$ is even. If $n$ is even, then $m$ and $m'$ are both even with
$\ord_{2} \left( m' \right) >\ord_{2}(m)>0$.

If $r=0$, then $k$ is odd and $2k+r$ is even. If $n$ is odd, then $m$ is
even and $m'$ is odd. If $n$ is even, then $m$ and $m'$ are both even with
$\ord_{2} \left( m \right) >\ord_{2}(m')>0$.

In each case, by Lemma~\ref{lem:2.8}(a) (which is applicable for $n$ odd since
$x(P)$ is assumed to be a rational square in this case), we have
$$
\left( A_{m}B_{m'}-A_{m'}B_{m} \right)^{2}
\leq B_{m+m'}B_{m-m'}.
$$

Taking the logarithm of both sides gives
\begin{equation}
\label{eq:2110}
2\log \left| A_{m}B_{m'}-A_{m'}B_{m} \right|
\leq \log \left( B_{n} \right) +\log \left( B_{m-m'} \right).
\end{equation}

Assume that the term $B_{n}$ does not have a primitive divisor. Then, by
Lemma~\ref{lem:2.7}, we have
\begin{equation}
\label{eq:2120}
\log \left( B_{n} \right) \leq \eta(n)+2n^{2}\rho(n) \widehat{h}(P)+\omega(n)K.
\end{equation}

Lemma~\ref{lem:2.6} gives
\begin{align}
\log \left( B_{m-m'} \right) & \leq h((m-m')P) \notag \\
&\le 2\widehat{h}((m-m')P)+K = 2(m-m')^{2}\widehat{h}(P)+K.  \label{eq:2130}
\end{align}

Combining \eqref{eq:2120} and \eqref{eq:2130} with \eqref{eq:2110} gives
\begin{align}
&2\log|A_{m}B_{m'}-A_{m'}B_{m}|  \notag\\
&\le \eta(n)+2n^{2}\rho(n) \widehat{h}(P)+\omega(n)K+2(m-m')^{2}\widehat{h}(P)+K.\label{eq:2140} 
\end{align}

Lemma~6.1(a) of \cite{VY1} implies that $A_{m}$ and $A_{m'}$ are both squares,
so we can write $A_{m}=a_{m}^{2}$, $A_{m'}=a_{m'}^{2}$
$B_{m}=b_{m}^{2}$ and $B_{m'}=b_{m'}^{2}$. Thus
\begin{eqnarray*}
2\log \left| A_{m}B_{m'}-A_{m'}B_{m} \right|
&   =  & 2\log \left| a_{m}^{2}b_{m'}^{2}-a_{m'}^{2}b_{m}^{2} \right| \\
& \geq & 2\log \left( \left| a_{m}b_{m'}|+|a_{m'}b_{m} \right| \right) \\
& \geq & 2\log \left( |a_{m}|+|b_{m}| \right) \\
& \geq & 2 \log \max \{|a_{m}|, |b_{m}|\} \\
&   =  & h(mP) \geq 2\widehat{h}(mP)-L,
\end{eqnarray*}
recalling from Lemma~\ref{lem:2.8} that $A_{m}B_{m'}-A_{m'}B_{m} \neq 0$. Note
that the last inequality is obtained by Lemma~\ref{lem:2.6} and the definition
of $L$ in (\ref{eq:L-value}). Since $\widehat{h}(mP)=m^{2}\widehat{h}(P)$, we
have
$$
2\log \left| A_{m}B_{m'}-A_{m'}B_{m} \right|
\geq 2m^{2}\widehat{h}(P)-L.
$$

Combining this estimate and \eqref{eq:2140} gives
\begin{eqnarray*}
&      & 2m^{2}\widehat{h}(P)-L \\
& \leq & \eta(n)+2n^{2}\rho(n)\widehat{h}(P)+\omega(n)K
+2(m-m')^{2}\widehat{h}(P)+K.
\end{eqnarray*}

Substituting $m=2^{e}(2N+r)/3$ and $m'=2^{e}(N-r)/3$ gives
\begin{align}
\eta(n)+\omega(n)K+K+L
& \geq 2\left(\frac{1}{3}-\frac{r^{2}}{3N^{2}}-\rho(n)\right)\widehat{h}(P)n^{2}  \notag \\
& \geq 2\left(\frac{1}{3}-\frac{1}{3N^{2}}-\rho(n)\right)\widehat{h}(P)n^{2}.\label{eq:2150}  
\end{align}

\subsubsection{Proof of part~(b)}

Assume that $n$ is a positive integer divisible by an odd prime $p$ and $x(P)$
is a rational square, if $n$ is odd. Write $n=pk$ for some positive integer $k$.
Assume that $B_{n}$ does not have a primitive divisor. Then by Lemmas~\ref{lem:2.7}
and \ref{lem:2.8}(a) (observe that $\ord_{2}(p-1) \neq \ord_{2}(p+1)$ so
\ref{lem:2.8}(a) is also applicable if $n$ is even), we have
\begin{align}
& 2 \log \left| A_{(p+1)k/2}B_{(p-1)k/2}-A_{(p-1)k/2}B_{(p+1)k/2} \right| \notag \\
& \leq \log \left( B_{n} \right) + \log \left( B_{k} \right)
\leq \eta(n)+2n^{2}\rho(n)\widehat{h}(P)+\omega(n)K +\log \left( B_{k} \right). \label{eq:2170}
\end{align}

On the other hand, 
\begin{eqnarray}
&      & 2 \log \left| A_{(p+1)k/2}B_{(p-1)k/2}-A_{(p-1)k/2}B_{(p+1)k/2} \right| \notag \\
&  =   & 2\log \left| a_{(p+1)k/2}^{2}b_{(p-1)k/2}^{2}-a_{(p-1)k/2}^{2}b_{(p+1)k/2}^{2} \right| \notag \\
&  =   & 2\log \left| \left| a_{(p+1)k/2}b_{(p-1)k/2} \right|
                        - \left| a_{(p-1)k/2}b_{(p+1)k/2} \right| \right| \notag \\
&      & + 2\log \left( \left| a_{(p+1)k/2}b_{(p-1)k/2} \right|
                          + \left| a_{(p-1)k/2}b_{(p+1)k/2} \right| \right) \notag\\
& \geq & 2\log \left| b_{k} \right|
         + 2\log \left( \left| a_{(p+1)k/2} \right| + \left| b_{(p+1)k/2} \right| \right) \qquad \text{since $b_{k} \mid b_{(p \pm 1)k/2}$,} \notag\\
& \geq & \log \left( B_{k} \right) + h(((p+1)k/2)P) \notag \\
& \geq & \log \left( B_{k} \right) + 2((p+1)k/2)^{2}\widehat{h}(P)-L.  \label{eq:2180}
\end{eqnarray}

Combining \eqref{eq:2170} and \eqref{eq:2180} gives
$$
2((p+1)k/2)^{2}\widehat{h}(P) - L
\leq \eta(n)+2n^{2}\rho(n)\widehat{h}(P) + \omega(n)K.
$$

Substituting $k=n/p$, we obtain
\begin{equation}
\label{eq:2190}
2 \left( \frac{(p+1)^{2}}{4p^{2}}-\rho(n) \right)\widehat{h}(P)n^{2}
\leq \eta(n) + \omega(n)K + L.
\end{equation}

We have thus completed the proof of part~(b).
\hfill $\square$


\subsection{Proof of Corollary~\ref{cor:1.4}}

To prove Corollary~\ref{cor:1.4}, we use Robin's estimate for $\omega(n)$
(see Th\'eor\`eme 11 of \cite{Robin}):
\begin{equation}
\label{eq:2210}
\omega(n) < \frac{1.3841\log (n)}{\log \log (n)}
\qquad \text{for all $n \ge 3$}.
\end{equation}

Furthermore, we use the following estimate for $\rho(n)$: 
\begin{eqnarray}
\label{eq:rho-bnd}
\rho(n) & < & \sum_{p < 10^{6}} p^{-2} + \left( \zeta(2) - \sum_{m \leq 10^{6}} m^{-2} \right)
<0.452248 + 0.000001 \notag \\
& < & 0.45225,
\end{eqnarray}
where the first sum is over primes, $p$, and the second sum over positive
integers, $m$.

\subsubsection{Proof of Corollary~\ref{cor:1.4}(a)}

Let $P \in E_{a}(\bQ)$ be a point of infinite order. Let $n \geq 3$ be an odd
integer, and assume that $x(P)$ is a rational square. We will distinguish
three cases.

{\it Case 1.} Assume that $n$ is not divisible by either $3$ or $5$. Then
$n \geq 7$ and $\rho(n) < 0.45225-1/4-1/9-1/25 < 0.052$. Here we apply
Proposition~\ref{prop:1.3}(a), so we have $N=n$ and
$$
2 \left( \frac{1}{3}-\frac{1}{3N^{2}}-\rho(n) \right) > 0.549,
$$
and the Corollary follows in this case.

{\it Case 2.} Assume that $n$ is divisible by $3$. Then $\rho(n) < 0.45225-1/4 < 0.2023$.
Here we apply Proposition~\ref{prop:1.3}(b) with $p=3$, so
$$
2 \left( \frac{(p+1)^{2}}{4p^{2}} - \rho(n) \right)
= 2 \left( \frac{4}{9} - \rho(n) \right)
> 0.484,
$$
and the Corollary follows in this case.

{\it Case 3.} Assume that $n$ is divisible by $5$, but not by $3$. Then
$\rho(n) < 0.45225-1/4-1/9<0.092$. Here we apply Proposition~\ref{prop:1.3}(b)
with $p=5$, so
$$
2 \left( \frac{(p+1)^{2}}{4p^{2}} - \rho(n) \right)
= 2 \left( \frac{9}{25} - \rho(n) \right)
> 0.536.
$$

Therefore, the Corollary follows in this case too, completing the proof of
part~(a).


\subsubsection{Proof of Corollary~\ref{cor:1.4}(b)}

Let $n$ be a positive even integer and assume that $B_{n}(P)$ does not have a
primitive divisor. If $n$ is a power of two, then $n \le 4$, so by excluding
these values of $n$ in the hypotheses of the Corollary, we may assume here that
$n$ is not a power of two. We will distinguish two cases.

{\it Case 1.} Assume that $n$ is not divisible by $3$ and $5$. From
(\ref{eq:rho-bnd}), $\rho(n) < 0.45225-1/9-1/25 < 0.302$.

Here we apply Proposition~\ref{prop:1.3}(a) and write $n=2^{e}N$ with $e \geq 1$
and $N \geq 7$ odd. In this way, we obtain
$$
2\left( \frac{1}{3}-\frac{1}{3N^{2}}-\rho(n) \right) > 0.049,
$$
and the Corollary follows in this case.

{\it Case 2.} Assume that $n$ is divisible by $3$ and not by $5$. Then
$\rho(n) < 0.41225$. Here we apply Proposition~\ref{prop:1.3}(b) with $p=3$,
so
$$
2 \left( \frac{(p+1)^{2}}{4p^{2}} - \rho(n) \right)
= 2 \left( \frac{4}{9} - \rho(n) \right)
> 0.064,
$$
and the Corollary follows in this case.

We have thus completed the proof. 
\hfill $\square$

\section{Proof of Theorem~\ref{thm:1.6}}

\subsection{$n$ divisible by $5$}

\begin{lemma}
\label{lem:7}
Let $P \in E_{a}(\bQ)$ be a point of infinite order, and
$\left( B_{n} \right)=\left( B_{n}(E_{a},P) \right)$ an elliptic divisibility
sequence. For all positive integers $m$, $B_{5m}$ has a primitive
divisor.
\end{lemma}

\begin{remark}
This technique can be applied for other values than $5$, provided that
$\psi_{n}(a,x)$ is reducible over $\bQ$ (e.g., $n=13$, $17$,\ldots, but
not $n=7$, $11$, $19$, \ldots, for which $\psi_{n}(a,x)$ are irreducible).
\end{remark}

\begin{proof}
We can handle $m=1,2,3,4$ and $5$ using the arguments in Section~3 of
\cite{Ingram} as they readily generalise to our curves.

For $m \geq 6$, we follow the idea of Ingram \cite[Lemma~7]{Ingram}.

For a point $P=(x,y) \in E_{a}(\bQ)$, we can factor $\psi_{5}(P)=\psi_{5}(a,x)$,
the 5-th division polynomial, as
$$
\left( a^{2}+2ax^{2}+5x^{4} \right)
\left( x^{8}+12ax^{6}-26a^{2}x^{4}-52a^{3}x^{2}+a^{4} \right).
$$

Writing $x([m]P)=u/v^{2}$ with $\gcd(u,v)=1$, then
$$
x([5m]P) = \frac{\phi_{5}(a, u/v^{2})}{\psi_{5}^{2}(a, u/v^{2})},
$$
where $\psi_{5}(a,x)$ is a binary form in $a$ and $x^{2}$ of degree $12$, while
$\phi_{5}(a,x)$ is $x$ times a binary form in $a$ and $x^{2}$ of degree $24$.
Hence $v^{50}\phi_{5}(a, u/v^{2}), v^{50}\psi_{5}^{2}(a, u/v^{2}) \in \bZ$ and
$$
B_{5m} = v^{50}\psi_{5}^{2}(a,u/v^{2})/g,
$$
where $g=\gcd \left( v^{50}\phi_{5}(a, u/v^{2}), v^{50}\psi_{5}^{2}(a, u/v^{2}) \right)$.

We can write $v^{25}\psi_{5}(a,u/v^{2})=
v f_{5,1} \left( u^{2}, v^{4} \right)
f_{5,2} \left( u^{2}, v^{4} \right)$,
where
$$
f_{5,1}(x,y)= x^{4}+12ax^{3}y-26a^2x^{2}y^{2}-52a^{3}xy^{3}+a^{4}y^{4}
$$
and
$$
f_{5,2}(x,y)= 5x^{2}+2axy+a^{2}y^{2}.
$$

\underline{Lower bounds in terms of $y$.}
The roots of $f_{5,1}(x,1)$ are at $r_{1,1}=-13.6275\ldots a$,
$r_{1,2}=-1.3167\ldots a$, $r_{1,3}=0.0190\ldots a$ and $r_{1,4}=2.9252\ldots a$.

The roots of $f_{5,2}(x,1)$ are $r_{2,1}=(-1-2i)a/5$ and $r_{2,2}=(-1+2i)a/5$.
Hence $f_{5,2}(x,y)=5\left( x-r_{2,1}y \right)\left( x-r_{2,2}y \right)$.

For any integers, $x$ and $y$, with $y \neq 0$, there will be an $i$ and a $j$
such that $|x/y-r_{i,j}|$ is minimal.

Suppose the closest root to $z=x/y$ is either $r_{2,1}$ or $r_{2,2}$. The
nearest roots of $f_{5,1}(x,1)$ to $r_{2,1}$ and $r_{2,2}$ are $r_{1,2}$ and
$r_{1,3}$, so by solving $(z-(-1.3167a))^{2}=(z-(-a/5))^{2}+(2a/5)^2$ and
$(z-(0.019a))^{2}=(z-(-a/5))^2+(2a/5)^{2}$ respectively for $z$, we find
that $x/y$ must lie between $-0.455a$ and $-0.687a$.

For such $x$ and $y$,
$$
\left| f_{5,1}(x,y) \right|
> y^{4} \left| -0.687a-r_{1,1} \right|
\left| -0.687a-r_{1,2} \right|
\left| -0.455a-r_{1,3} \right|
\left| -0.455a-r_{1,4} \right|
> 13a^{4}y^{4}
$$
and
$$
\left| f_{5,2}(x,y) \right| > 5y^{2} \left| -0.455a-r_{2,1} \right|
\left| -0.455a-r_{2,2} \right|
1.125a^{2}y^{2}.
$$

In this case,
$$
\left| y^{1/4} \right| \left| f_{5,1} \left( x, y \right) \right|
\left| f_{5,2} \left( x, y \right) \right|
> 14.6a^{6}|y|^{25/4}.
$$

If the closest root to $x/y$ is one of the $r_{1,j}$'s, then
$\left|f_{5,2}(x,y) \right| \geq 5 I_{y}^{2}=4a^{2}y^{2}/5$ where
$I_{y}=|2ay/5|$ is the absolute value of the imaginary part of $r_{2,1}$
and $r_{2,2}$.

In this case,
\begin{equation}
\label{eq:v-lb}
\left| y^{1/4} \right| \left| f_{5,1} \left( x, y \right) \right|
\left| f_{5,2} \left( x, y \right) \right|
\geq (4/5)a^{2}|y|^{9/4}
> 3.8|y|^{9/4}/|a|^{1/4},
\end{equation}
since $|a| \geq 2$ and $\left| f_{5,1} \left( x,y \right) \right| \geq 1$ here.

Hence this last inequality always holds (including for $y=0$).

\underline{Lower bounds in terms of $x$.}
The roots of $f_{5,1}(1,y)$ are at $s_{1,1}=-0.7594\ldots/a$,
$s_{1,2}=-0.0733\ldots/a$, $s_{1,3}=0.3418\ldots/a$ and $s_{1,4}=52.4909\ldots/a$.

And the roots of $f_{5,2}(1,y)$ are $s_{2,1}=(-1-2i)/a$ and $s_{2,2}=(-1+2i)/a$.
Hence $f_{5,2}(x,y)=a^{2}\left( y-s_{2,1}x \right)\left( y-s_{2,2}x \right)$.

For any integers, $x$ and $y$, with $x \neq 0$, there will be an $i$ and a $j$
such that $|y/x-s_{i,j}|$ is minimal.

Suppose the closest root to $y/x$ is either $s_{2,1}$ or $s_{2,2}$. The nearest
root of $f_{5,1}(1,y)$ is $s_{1,1}$. By solving
$(z-(-0.7594/a))^{2}=(z-(-1/a))^{2}+(2/a)^2$ for $z=y/x$, we find
that $y/x$ must lie beyond $-9.192/a$.

For such $x$ and $y$, $|f_{5,2}(x,y)| > 71.1x^{2}$ and
$|f_{5,1}(x,y)| > 45,200x^{4}$.

In this case,
$$
\left| y^{1/4} \right| \left| f_{5,1} \left( x, y \right) \right|
\left| f_{5,2} \left( x, y \right) \right|
> 3.2 \cdot 10^{6} x^{6}\left| y^{1/4} \right|.
$$

Since $|y/x|>9.192/|a|$,
$$
\left| y^{1/4} \right| \left| f_{5,1} \left( x, y \right) \right|
\left| f_{5,2} \left( x, y \right) \right|
> 5.5 \cdot 10^{6} |x|^{25/4}/|a|^{1/4}.
$$

If the closest root to $y/x$ is one of the $s_{1,j}$'s, then
$\left| f_{5,2}(x,y) \right| \geq a^{2} (2/a)^{2}x^{2}=4x^{2}$.

If $|y/x|>0.045/|a|$,
\begin{equation}
\label{eq:u-lb}
\left| y^{1/4} \right| \left| f_{5,1} \left( x, y \right) \right|
\left| f_{5,2} \left( x, y \right) \right|
> (0.045)^{1/4}|x|^{1/4}/|a|^{1/4} \cdot 4x^{2}
> 1.83 |x|^{9/4}/|a|^{1/4},
\end{equation}
since $\left| f_{5,1} \left( x,y \right) \right| \geq 1$.

If $|y/x| \leq 0.045/|a|$, we find that 
$\left| f_{5,1}(x,y) \right| > 0.3145x^{4}>0.374x^{4}/|a|^{1/4}$ and
$\left| f_{5,2}(x,y) \right| > 4.912x^{2}$, so
$$
\left| y^{1/4} \right| \left| f_{5,1} \left( x, y \right) \right|
\left| f_{5,2} \left( x, y \right) \right|
> 0.374 x^{4} \cdot 4.912x^{2}/|a|^{1/4}
\geq 1.83x^{6}/|a|^{1/4},
$$
for $|y| \geq 1$ and $|a| \geq 2$.

So \eqref{eq:u-lb} always holds.

Combining \eqref{eq:v-lb} and \eqref{eq:u-lb} with $x=u^{2}$ and $y=v^{4}$,
we find that
\begin{eqnarray}
\label{eq:poly-lb}
|v| \left| f_{5,1} \left( u^{2}, v^{4} \right) \right|
\left| f_{5,2} \left( u^{2}, v^{4} \right) \right|
&   >  & \min \left\{ 1.83 |u|^{4.5}/|a|^{1/4}, 3.8 |v|^{9}/|a|^{1/4} \right\} \nonumber \\
& \geq & 1.83|a|^{-1/4} \min \left\{ |u|^{4.5}, |v|^{9} \right\}.
\end{eqnarray}

Hence
$$
\log \left| B_{5m} \right| > 2\log(1.83) + 9h(u/v^{2})-(1/2)\log|a| -\log(g).
$$

Next, using the Maple command \verb+gcdex(psi[5]^2,phi[5],x,'s','t');+
and factoring the common denominator, we find that this gcd must divide
$2^{45}a^{24}$. As in the proof of Lemma~4 of \cite{Ingram}, we have
$g|2^{45}a^{24}$ and
$$
\log \left| B_{5m} \right| > 9h(u/v^{2})-(49/2)\log|a| - 29.983.
$$

From Lemma~\ref{lem:2.6},
$$
h(mP)=h(u/v^{2}) \geq 2\widehat{h}(mP)-(1/2)\log|a|-0.32,
$$
so
\begin{eqnarray*}
\log \left| B_{5m} \right|
& > & 18\widehat{h}(mP)-(9/2)\log|a|-2.88-(49/2)\log|a| - 29.983 \\
& > & 18m^{2}\widehat{h}(P)-29\log|a|-32.863,
\end{eqnarray*}
since $\widehat{h}(mP)=m^{2}\widehat{h}(P)$.

From Lemma~\ref{lem:2.7} applied with $n=5$ and using $mP$ here in place
of $P$ there:
$$
\log \left| B_{5m} \right| \leq 2\log(5)+2\widehat{h}(mP)+(1/2)\log|a| + 0.52
< 2m^{2}\widehat{h}(P)+(1/2)\log|a| + 3.739.
$$

Hence
\begin{equation}
\label{eq:test}
16m^{2}\widehat{h}(P) \leq (59/2)\log|a|+36.602.
\end{equation}

Unless $a \equiv 4 \bmod 16$, from Lemma~\ref{lem:lang},
\begin{equation}
\label{eq:5-lang}
\widehat{h}(P) \geq (1/16) \log |a| + (1/4)\log(2).
\end{equation}

Furthermore, if $a \equiv 4 \bmod 16$ and $x(2P) \not\in \bZ$, then
$\delta \geq 2$ if $\ord_{2}(x(2P))=0$ and $\delta \geq 3$ otherwise. In
both cases, we find from Proposition~1.1 of \cite{VY1} that \eqref{eq:5-lang}
holds.

Hence excluding the case of $a \equiv 4 \bmod 16$ and $x(2P) \in \bZ$,
we have
$$
m^{2} \left( \log |a| + 4\log(2) \right) \leq (59/2)\log|a|+36.602.
$$

For $m \geq 6$, this can never hold.

For $a \equiv 4 \bmod 16$ and $a>0$, we have
$$
m^{2} \left( \log |a| - 2\log(2) \right) \leq (59/2)\log|a|+36.602,
$$
or
$$
|a|
< \exp \left( \left( 36.602+2m^{2}\log(2) \right)/\left( m^{2}-(59/2) \right) \right).
$$

For $m \geq 6$, this last inequality is false once $a>\exp(13.31)$ and hence
for $a>604,200$, since the right-hand side is monotonically decreasing for
$m \geq 6$.

For $a \equiv 4 \bmod 16$ and $a<0$, we have
$$
m^{2} \left( \log |a| - \log(2) \right) \leq (59/2)\log|a|+36.602,
$$
or
$$
|a|
< \exp \left( \left( 36.602+m^{2}\log(2) \right)/\left( m^{2}-(59/2) \right) \right).
$$

For $m \geq 6$, this last inequality is false once $|a|>\exp(9.471)$ and hence
for $|a|>13,000$, since the right-hand side is monotonically decreasing for $m \geq 6$.

For $m \geq 7$, the required inequalities hold for $a<-37$ and $a>213$.

\underline{Search.}
For each of the remaining values of $a$ with $a \equiv 4 \bmod 16$, we use
PARI to search for possible counterexamples to our lemma.

We search for any points, $P$, with $x(2P) \in \bZ$, satisfying \eqref{eq:test}
for $m=6$ and then check to ensure that $B_{30}(P)$ always has a primitive
divisor for such points.

Using Lemma~\ref{lem:2.6} and since $\widehat{h}(2P)=4\widehat{h}(P)$,
we require
\begin{eqnarray*}
h(2P) & = & \log|x(2P)| \leq 2 \left( (1/4) \log|a| + 0.26 + 4\frac{29.5\log|a|+36.602}{16 \cdot 36} \right) \\
& < & 0.91\log|a| + 1.03.
\end{eqnarray*}

Furthermore, $x(2P)$ must be a square by Lemma~6.1(a) of \cite{VY1}.

In this way, we were able to check all the remaining values in 4 minutes using
PARI on an ordinary laptop. We found 29 pairs $(a,P)$ such that \eqref{eq:test}
holds and $a$ is fourth-power-free (6 with $-9996 \leq a \leq -12$ and $23$
with $180 \leq a \leq 515,508$).

For all of them, $B_{30}(P)$ has a primitive divisor. This completes the proof
for $m=6$.

Since the left-hand side of \eqref{eq:test} is monotonically increasing in $m$,
any pair $(a,P)$ satisfying \eqref{eq:test} for $m>6$, must be among these 29
examples. As we noted above, we must have $-37 \leq a \leq 213$ if $B_{5m}(P)$
does not have a primitive divisor for $m>6$, so this leaves only $a=-12,180$.

For $a=-12$, the point is $P=(-2,4)$ and $\widehat{h}(P)=0.1252\ldots$.
For $m \geq 8$, \eqref{eq:test} is no longer satisfied in this case and we check
that $B_{35}(P)$ has a primitive divisor.

For $a=180$, the point is $P=(-6,36)$ with $\widehat{h}(P)=0.2564\ldots$.
We proceed in the same way for this case.
\end{proof}

\subsection{$n$ odd, $n>7$}

Let $n$ be a positive odd integer and assume that either $x(P)$ is a rational
square or $x(P)<0$. Suppose that the term $B_{n}(P)$ does not have a primitive
divisor.

\subsubsection{$x(P)<0$}

This only occurs for $a<0$. In this case, $E_{a}(\bR)$ has two components and we
are considering points, $P$, on the non-identity component of $E_{a}(\bR)$.

\begin{remark}
Since Lemma~\ref{lem:2.7} applies for any elliptic curve and since we can
obtain results like Lemma~\ref{lem:2.6} for any elliptic curve (see, for
example, Proposition~5.18(a) and Theorem~5.35(c) of \cite{S-Z}), if $E$ is
an elliptic curve such that $E(\bR)$ has two components and if we have
an explicit version of Lang's conjecture for $E(\bQ)$, then the same
idea can be applied to bound from above the odd indices $n$, such that
$B_{n}$ has no primitive divisor when $P$ is on the non-identity component of
$E(\bR)$.
\end{remark}

If $B_{n} \geq \left| A_{n} \right|$, then
\begin{equation}
\label{eq:bn-lbnd1}
\log \left( B_{n} \right) = h(x(nP)) \geq 2\widehat{h}(nP)- \frac{1}{2} \log |a| - 0.32
= 2n^{2}\widehat{h}(P)- \frac{1}{2} \log |a| - 0.32,
\end{equation}
from Lemma~\ref{lem:2.6}.

If $B_{n} < \left| A_{n} \right|$, then
$$
-\sqrt{|a|} \leq x(nP) < 0,
$$
since $x(P)<0$ and $n$ odd implies $x(nP)<0$.

Thus $\left| A_{n}/B_{n} \right| \leq \sqrt{|a|}$, so
$$
\log \left| A_{n} \right| - \frac{1}{2} \log |a| \leq \log \left( B_{n} \right).
$$

Therefore, from Lemma~\ref{lem:2.6},
\begin{equation}
\label{eq:bn-lbnd2}
\log \left( B_{n} \right) > h(x(nP)) - \frac{1}{2} \log |a|
\geq 2n^{2}\widehat{h}(P)- \log |a| - 0.32.
\end{equation}

Since this lower bound is weaker than (\ref{eq:bn-lbnd1}), we shall use it in
what follows.

Applying Lemma~\ref{lem:2.7}, we have
\begin{equation}
\label{eq:n-odd-x-neg-bnd}
2n^{2} \left( 1-\rho(n) \right) \widehat{h}(P)
< \eta(n)+\omega(n)K + \log |a| + 0.32.
\end{equation}

From Lemma~\ref{lem:lang},
$$
\frac{n^{2}}{8} \left( 1-\rho(n) \right) \left( \log |a| - 2\log(2) \right)
< \eta(n)+\omega(n)K + \log |a| + 0.32.
$$

For $|a| \geq 8$,
\begin{eqnarray*}
\frac{1}{\log |a| - 2\log(2)} & < & 1.443, \quad
\frac{K}{\log |a| - 2\log(2)} < 2.251, \\
\frac{\log |a| + 0.32}{\log |a| - 2\log(2)} & < & 3.462.
\end{eqnarray*}

Hence,
\begin{equation}
\label{eq:n-odd-x-neg}
\frac{n^{2}}{8} \left( 1-\rho(n) \right)
< 1.443\eta(n)+2.251\omega(n) + 3.462.
\end{equation}

Since $n$ is odd, from \eqref{eq:rho-bnd}, $\rho(n)<0.20225$ and applying
\eqref{eq:2210}
$$
0.0997n^{2} < 2.886\log(n) + 3.116\frac{\log (n)}{\log \log (n)} + 3.462.
$$

Using this inequality, we obtain the bound $n<14.01$, so $n \leq 13$.
We can directly calculate both sides of (\ref{eq:n-odd-x-neg}) and find
that it fails for $n=9$ and $n=11$, proving our desired result for $|a| \geq 8$.

Using PARI, we find that for $-8 < a <0$, $E_{a}(\bQ)$ is non-trivial only for
$a=-2, -5, -6$ and $-7$. Each of these is of rank $1$ and the generator,
$P$, satisfies $\widehat{h}(P) \geq 0.3043\ldots$ (attained for $a=-2$).
Applying this lower bound for the height and $\rho(n)<0.20225$ to
(\ref{eq:n-odd-x-neg-bnd}), we find that
$$
0.485n^{2}
> 2\log(n)+2.16 \frac{\log(n)}{\log(\log(n))}+2.4
> 2\log(n)+\omega(n)K+\log|a|+0.32
$$
for $-8 \leq a \leq -2$ and $n>5.2$, as required.

\noindent
{\bf Note:} recall that the height function, ellheight, in PARI returns a value that is
twice our height here.

\subsubsection{$x(P)$, a rational square}

From Lemma~\ref{lem:lang}, we have
\begin{equation}
\label{eq:350}
\widehat{h}(P) \geq \frac{1}{16}\log |a| - \frac{1}{8}\log(2).
\end{equation}

Assume further that $|a| \geq 33$. Then
\begin{eqnarray}
\label{eq:360}
\frac{1}{\log |a|-2\log(2)} & < & 0.474, \hspace{1.0mm}
\frac{K}{\log |a|-2\log(2)} < 1.075, \nonumber \\
\frac{L}{\log |a|-2\log(2)} & < & 0.981.
\end{eqnarray}

Substituting \eqref{eq:350} into Corollary~\ref{cor:1.4}(a) shows that
$$
0.484 \left( \frac{1}{16} \log |a| - \frac{1}{8}\log(2) \right) n^{2} \\
< 2 \log (n) + \frac{1.3841 \log (n)}{\log \log (n)}K + K + L
$$
must hold if $B_{n}(P)$ does not have a primitive divisor.

Dividing both sides of this equation by $\log |a|-2\log(2)$ and
substituting the estimates \eqref{eq:360} yields
\begin{eqnarray*}
0.030 n^{2} & < & 0.948\log (n) + 1.075\frac{1.3841 \log (n)}{\log \log (n)} + 2.056 \\
& < & \log (n) \left( 0.948 + \frac{1.488}{\log \log (n)} \right) + 2.056.
\end{eqnarray*}

Using this inequality, we obtain the bound $n<17.13$, so $n \leq 17$,
if $B_{n}(P)$ does not have a primitive divisor.

We will next give the better bounds by using the inequalities of
Proposition~\ref{prop:1.3}(b).

If $n$ is odd and divisible by $p$, then, by Proposition~\ref{prop:1.3}(b),
(\ref{eq:350}) and (\ref{eq:360}),
\begin{eqnarray*}
0 & < & \frac{1}{8}\left( \frac{(p+1)^{2}}{4p^{2}}-\rho(n) \right) n^{2}
\leq \frac{\eta(n)+\omega(n)K+L}{\log |a|-2\log(2)} \\
& < & 0.474 \eta(n) + 1.075 \omega (n) + 0.981.
\end{eqnarray*}

Using this inequality, we can eliminate $n=9$, $11$, $13$, $15$ and $17$ (with
$p=3$, $11$, $13$, $3$ and $17$, respectively) for $|a| \geq 33$.

Next using Lemma~\ref{lem:2.6}, if $x(P)$ is a rational square and
$\widehat{h}(P) \leq (1/4)\log|a|$, then
$h(P) \leq (1/2)\log|a| + (1/2)\log|a|+0.52 = \log|a|+0.52$. That is, writing
$x(P)=r/s$, $\max \left( |r|, |s| \right) < 1.7 |a|$. Using this bound on
$x(P)$ we can enumerate all points with $x(P)$ a rational square,
$\widehat{h}(P) \leq (1/4)\log|a|$ and $|a| \leq 32$ (in PARI, say):
\begin{table}[h]
\begin{tabular}{|r c c | r c c |}
\hline
  $a$ & $P$           & $\widehat{h}(P)$ \\ \hline
$-12$ & $(4, \pm 4)$  & $0.5011\ldots$   \\
  $3$ & $(1, \pm 2)$  & $0.2505\ldots$   \\
 $15$ & $(1, \pm 4)$  & $0.5673\ldots$   \\
 $20$ & $(4, \pm 12)$ & $0.6355\ldots$   \\
\hline
\end{tabular}
\caption{}
\label{table:n-odd}
\end{table}

Repeating the above arguments with $|a| \geq 2$ (noting the $E_{-1}(\bQ)$ and
$E_{1}(\bQ)$ contain only torsion points) and $\widehat{h}(P) \geq (1/4) \log |a|$,
we obtain
$$
0.12 n^{2} < \log (n) \left( 2.886 + \frac{1.732}{\log \log (n)} \right)
+ 2.213,
$$
from Corollary~\ref{cor:1.4}(a),
and
$$
0 < \frac{1}{2} \left( \frac{(p+1)^2}{4p^{2}} - \rho(n) \right) n^{2}
< 1.443\eta(n) + 1.251\omega(n) + 0.962,
$$
from Proposition~\ref{prop:1.3}(b).

From the first of these inequalities, we find that $n<10.8$ and so $n \leq 9$
for all $a$ and $P \in E_{a}(\bQ)$ with $|a| \geq 2$ and $\widehat{h}(P) \geq (1/4)\log |a|$.
The second inequality allows us to eliminate $n=9$. Hence we find that
$n \leq 7$ for all elliptic divisibility sequences such that $B_{n}(P)$ does
not have a primitive divisor, with the exception of those generated by the $8$
points, $P$, in Table~\ref{table:n-odd}.

Substituting $a=3$ and $\widehat{h}(P)=0.2505\ldots$ into the inequality in
Corollary~\ref{cor:1.4}(a), we find that
$$
0.484n^{2} < \left( 7.984 + \frac{5.912}{\log \log(n)} \right) \log(n)+7.744.
$$

Using this inequality, we find that $n \leq 9$. For $n=9$, we use
Proposition~\ref{prop:1.3}(b) with $p=3$. The left-hand side exceeds the
right-hand side, and so we can eliminate $n=9$ too.

We proceed in the same way for $a=-12$ and $15$ and $20$.

This completes the proof for $n$ odd.

\subsection{$n$ even, $n \geq 20$}

Let $n$ be a positive even integer and not a power of two. Assume that $B_{n}(P)$
does not have a primitive divisor.

We suppose that $|a| \geq 384$. Then
\begin{eqnarray}
\label{eq:370}
\frac{1}{\log |a|-2\log(2)} & < & 0.22, \hspace{1.0mm}
\frac{K}{\log |a|-2\log(2)} < 0.766, \nonumber \\
\frac{L}{\log |a|-2\log(2)} & < & 0.722.
\end{eqnarray}

From \eqref{eq:350}, \eqref{eq:370} and
using the same argument as for $n$ odd, with Corollary~\ref{cor:1.4}(b), we have
$$
0.003 n^{2} < \left( 0.44 + \frac{1.061}{\log \log (n)} \right) \log (n) + 0.722.
$$

By using this inequality, we obtain the bound $n < 42.4$, so $n \leq 42$.

From Proposition~\ref{prop:1.3}(a), \eqref{eq:350} and \eqref{eq:370}, we have
$$
\frac{1}{8} \left( \frac{1}{3} - \frac{1}{3N^{2}} - \rho(n)\right) n^{2}
< 0.22\eta(n) + 0.766\omega(n) + 1.488.
$$

Using this inequality, we obtain $n \leq 20$ or $n=24$, $30$, $32$, $36$, $42$,
for $|a| \geq 384$. Furthermore, recall that we need not consider $n=30$, since
we have already eliminated those $n$ divisible by $5$.

Applying Proposition~\ref{prop:1.3}(b), \eqref{eq:350} and \eqref{eq:370}, we have
$$
\frac{1}{8} \left( \frac{(p+1)^{2}}{4p^{2}} - \rho(n)\right) n^{2}
< 0.22\eta(n) + 0.766\omega(n) + 0.722.
$$

Using this inequality for the remaining values of $n$, we obtain $n < 20$
for $|a| \geq 384$ ($n \leq 14$, in fact).

Now assume that $-384 < a < 384$. Again, we can proceed as in the case of
$n$ odd. First, we use Corollary~\ref{cor:1.4}(b), Proposition~\ref{prop:1.3}(a)
and Proposition~\ref{prop:1.3}(b) to shows that if $B_{n}(P)$ has no primitive
divisor, then $n \leq 20$ for $|a| \geq 2$ and $\widehat{h}(P) \geq (1/4) \log |a|$.

Then we use Lemma~\ref{lem:2.6} to find all points on $E_{a}(\bQ)$ with
$|a| \leq 384$ and $0<\widehat{h}(P) \leq (1/4)\log|a|$ and calculate $B_{n}(P)$
until the inequality in Corollary~\ref{cor:1.4}(b) is no longer satisfied.
Using PARI, 410 such points, $P$, were found. Typically, we needed to check
$B_{n}(P)$ for $n$ up to $20-25$, although for $a=-12$, the search had to
continue to $n=52$. The entire calculation took 11 minutes on an ordinary laptop.

\subsection{$n=3$}

As mentioned in the proof of Lemma~\ref{lem:2.8}, from the arguments on pages~92--93
of \cite{Silv-Tate}, we can write
$P= \left( b_{1}M^{2}/e^{2}, b_{1}MN/e^{3} \right)$ in lowest terms. By the
duplication formula, we have
$$
x(2P)= \frac{\left( 2b_{1}M^{4}-N^{2} \right)^{2}}{4M^{2}N^{2}e^{2}}.
$$

\noindent
\underline{$x(P)$ is a square.} Applying Lemma~\ref{lem:2.8}(a) with $m=2$ and
$n=1$, we have
\begin{equation}
\label{eq:n3-a}
0 < \left( A_{2}B_{1} - A_{1}B_{2} \right)^{2} \leq B_{1}B_{3}.
\end{equation}

Since $B_{1}|B_{2}$, we can write $B_{2}=k^{2}B_{1}$ for some integer $k \geq 1$.
Substituting this expression into \eqref{eq:n3-a}, we obtain
\begin{equation}
\label{eq:n3-b}
0 < B_{1} \left( A_{2} - k^{2}A_{1} \right)^{2} \leq B_{3}.
\end{equation}

If $\left| A_{2} - k^{2}A_{1} \right| > 3$, then from \eqref{eq:n3-b}, this
implies that $B_{3}>3^{2}B_{1}$ holds. Therefore from Lemma~\ref{lem:2.3},
$B_{3}$ has a primitive divisor.

Assume that $\left| A_{2} - k^{2}A_{1} \right|
\leq 3$. Then writing $A_{1}=a_{1}^{2}$ and $A_{2}=a_{2}^{2}$ with
$a_{1}, a_{2} \geq 1$, we have
\begin{equation}
\label{eq:n3-c}
\left| A_{2} - k^{2}A_{1} \right|
= \left| \left( a_{2}-ka_{1} \right) \left( a_{2}+ka_{1} \right) \right| \leq 3.
\end{equation}

From Lemma~6.1(c) of \cite{VY1}, we have $\ord_{2} \left( B_{2} \right) \geq \ord_{2} \left( B_{1} \right) + 2$
and so $k \geq 2$. By the left-hand inequality of \eqref{eq:n3-b}, $a_{2} \neq ka_{1}$
and so $a_{2}+ka_{1} \leq 3$. Together, these two statements imply that
$a_{1}=a_{2}=1$ and $k=2$. Since $A_{1}=1$ and $\gcd \left( b_{1}, e \right)
= \gcd \left( M, e \right)=1$, it follows that $b_{1}=M=\pm 1$ and so
$P= \left( 1/e^{2}, \pm N/e^{3} \right)$. Therefore, 
$x(2P)= \left( \pm 2 - N^{2} \right)^{2}/ \left( 4N^{2}e^{2} \right)
=1/ \left( 4e^{2} \right)$, since $A_{2}=1$ and $B_{2}=k^{2}B_{1}=4e^{2}$.
Hence, $N= \pm 1, \pm 2$.

If $N=\pm 1$, then $P=\left( 1/e^{2}, \pm 1/e^{3} \right)$, which is
impossible. Next assume that $N=\pm 2$. Then $P= \left( 1/e^{2}, \pm 2/e^{3} \right)$.
Substituting $x=1/e^{2}$ and $y=\pm 2/e^{3}$ into $y^{2}=x^{3}+ax$, we obtain
$a=3$ and $e=\pm 1$. That is $P=(1, \pm 2)$. In this case, we observe that
$B_{2}=4$ and $B_{3}=9$. Thus $B_{3}$ has a primitive divisor.

\noindent
\underline{$x(P)<0$.} Applying Lemma~\ref{lem:2.8}(b) with $m=2$ and $n=1$,
we have
$$
0 < \left( A_{2}B_{1} - A_{1}B_{2} \right)^{2} \leq 4B_{1}B_{3}.
$$

Once again we can write $B_{2}=k^{2}B_{1}$ for some integer $k \geq 1$ and we
have
$$
0 < B_{1} \left( A_{2} - k^{2}A_{1} \right)^{2} \leq 4B_{3}.
$$

In fact, from Lemma~\ref{lem:2.8}(a), we have
$$
0 < B_{1} \left( A_{2} - k^{2}A_{1} \right)^{2} \leq B_{3},
$$
unless $a \equiv 4 \bmod 16$ and $\ord_{2}(x(P))=1$.

If $\left| A_{2} - k^{2}A_{1} \right| > 6$, then $B_{3}>3^{2}B_{1}$ holds.
Therefore from Lemma~\ref{lem:2.3}, $B_{3}$ has a primitive divisor.

Assume that $\left| A_{2} - k^{2}A_{1} \right| \leq 6$ with $2||A_{1}$ or
$\left| A_{2}-k^{2}A_{1} \right| \leq 3$ otherwise. Since $A_{1}<0$ and $A_{2}$
is a square, we have $\left( A_{1}, A_{2}, k \right)
= (-2,4,1)$, $(-2,1,1)$ or $(-1,1,1)$.

In the first case, it follows that $b_{1}=-2$ and $M=\pm 1$, so
$P= \left( -2/e^{2}, \pm 2N/e^{3} \right)$. Therefore, $x(2P)=\left( -4-N^{2} \right)^{2}/
\left( 4N^{2}e^{2} \right)=4/e^{2}$ since $A_{2}=4$ and $B_{2}=k^{2}B_{1}=e^{2}$.
Hence $(-4-N^{2})/(2N)=\pm 2$, that is $N=\pm 2$ and so $P= \left( -2/e^{2}, \pm 4/e^{3} \right)$.
In order for $P$ to be on the curve $E_{a}$, $a=-12/e^{4}$. This implies that
$e=\pm 1$ and $a=-12$. In this case, $B_{1}=B_{2}=1$ and $B_{3}=3$, so $B_{3}$
has a primitive divisor.

In the second case (i.e., $\left( A_{1}, A_{2}, k \right)=(-2,1,1)$), once again
we have $x(2P)=\left(-4-N^{2} \right)^{2}/ \left( 4N^{2}e^{2} \right)$ and
$B_{2}=k^{2}B_{1}=e^{2}$. Since $A_{2}=1$, $x(2P)=1/e^{2}$, so
$(-4-N^{2})/(2N)=\pm 1$. There are no such rational $N$ and hence this
case is impossible too.

Lastly, we consider $\left( A_{1}, A_{2}, k \right)=(-1,1,1)$. Here we have
$b_{1}=-1$ and $M=\pm 1$, so $P= \left( -1/e^{2}, \pm N/e^{3} \right)$. Therefore,
$x(2P)= \left( -2-N^{2} \right)^{2}/ \left( 4N^{2}e^{2} \right)=1/e^{2}$
since $A_{2}=1$ and $B_{2}=k^{2}B_{1}=e^{2}$. Hence $\left( -2-N^{2} \right)/(2N)
=\pm 1$. There are no such rational $N$ and hence this case is impossible too.


Thus $B_{3}$ always has a primitive divisor in this case too.

\subsection{$4 \leq n \leq 20$}
\label{subsect:ingram}

In Section~2 of \cite{Ingram}, Ingram showed that there are no solutions for
$n=5$ and $7$. In Section~3 of the same paper, he proves the same for $n=4$,
$6$, $10$, $12$, $14$, $18$ and $20$, provided that $a=-N^{2}$. However,
writing $x(P)=A/B^{2}$ and setting $X=A^{2}/(A^{2},a)$ and $Y=aB^{4}/(A^{2},a)$,
rather than Ingram's values, we find that the polynomials, $\Psi_{n}$, are still
reducible (although not always with as many factors as in the $a=-N^{2}$ case).
Hence his same arguments hold (basically, for each of the possible values of
the factors use a reduction technique to eliminate one of the variables and
then factor the resulting single-variable polynomial to solve for $X$).

One can easily show that if $B_{n}(E,P)$ has no primitive divisor, then neither
does $B_{\rad(n)}(E,(n/\rad(n))P)$.

This completes the proof.

Lastly, note that the results for these values of $n$ hold without any
conditions on $x(P)$ and thus provide evidence for our claim in Remark~\ref{rem:1.5}.


\bibliographystyle{amsplain}

\begin{thebibliography}{20}
\bibitem{Bilu}
Y. Bilu, G. Hanrot, P. Voutier (with an appendix by M. Mignotte),
\emph{Existence of primitive divisor of Lucas and Lehmer numbers},
J. reine angew. Math. 539 (2001), 75--122.


\bibitem{Carm}
R. D. Carmichael,
\emph{On the numerical factors of the arithmetic forms $\alpha^{n} \pm \beta^{n}$},
Ann. of Math. 15 (1913), 30--70.

\bibitem{Cremona}
J. E. Cremona, M. Prickett, S. Siksek,
\emph{Height difference bounds for elliptic curves over number fields},
J. Number Theory 116 (2006), 42--68.


\bibitem{Durst}
L. K. Durst,
\emph{Exceptional real Lehmer sequences},
Pacific J. Math. 9 (1959), 437--441.

\bibitem{Everest1}
G. Everest, G. Mclaren, T. Ward,
\emph{Primitive divisors of elliptic divisibility sequenses},
J. Number Theory 118 (2006), 71--89.

\bibitem{Everest2}
G. Everest, A. van der Poorten, I. Shparlinski, and T. Ward,
{\em Recurrence sequences}, volume 104 of {\em Mathematical Surveys and Monographs},
American Mathematical Society, Providence, RI, 2003.

\bibitem{Ingram}
P. Ingram,
\emph{Elliptic divisibility sequences over certain curves},
J. Number Theory 123 (2007), 473--486.

%
\bibitem{Robin}
G. Robin,
\emph{Estimation de la fonction de Tchebychef $\theta$ sur le $k$-i\`{e}me
nombre premier et grandes valeurs de la fonction $\omega(n)$ nombre de diviseurs
premiers de $n$},
Acta Arith. XLII (1983), 367--389.

\bibitem{S-Z}
S. Schmitt, H. G. Zimmer,
\emph{Elliptic Curves: A computational Approach},
de Gruyter Studies in Mathematics 31, Walter de Gruyter Inc, 2004 


\bibitem{Silv2}
J.~H. Silverman,
\emph{The Arithmetic of Elliptic Curves},
Graduate Texts in Math. 106, Springer-Verlag, New York, 1986.


\bibitem{Silv4}
J.~H. Silverman,
\emph{Wieferich's criterion and the abc-conjecture},
J. Number Theory 30 (1988), 226--237.

\bibitem{Silv5}
J.~H. Silverman,
\emph{Advanced Topics in the Arithmetic of Elliptic Curves},
Graduate Texts in Math. 151, Springer-Verlag, New York, 1994.

\bibitem{Silv-Tate}
J.~H. Silverman, J. Tate,
\textit{Rational Points on Elliptic curves},
Undergraduate Texts in Math., Springer-Verlag, New York, 1992.


\bibitem{VY1}
P. Voutier and M. Yabuta,
\emph{Lang's conjecture and sharp height estimates for the elliptic curves $y^{2}=x^{3}+ax$}, 
(submitted).

\bibitem{Ward1}
M. Ward,
\emph{Memoir on elliptic divisibility sequences}, 
Amer. J. Math. 70 (1948), 31--74.

\bibitem{Ward2}
M. Ward,
\emph{The intrinsic divisors of Lehmer numbers},
Ann. of Math. (2) 62 (1955), 230--236.

\bibitem{Yabuta}
M. Yabuta,
\emph{Primitive divisors of certain elliptic divisibility sequences},
Exp. Math. 18(3) (2009), 303--310.


\bibitem{Zsigmondy}
K. Zsigmondy,
\emph{Zur Theorie der Potenzreste},
Monatsh. Math. 3 (1892), 265--284.
\end{thebibliography}

\end{document}